\documentclass[]{article}
\usepackage[utf8]{inputenc}
\usepackage[T1]{fontenc}

\usepackage[fleqn]{amsmath}
\usepackage{verbatim,amsthm}
\usepackage{graphicx}
\usepackage[colorinlistoftodos]{todonotes}
\usepackage[colorlinks=true, allcolors=blue]{hyperref}
\usepackage{amssymb,xcolor,amscd}
\usepackage{float}
\usepackage{enumitem}
\usepackage{listings}
\usepackage{authblk}
\usepackage{dsfont}

\newcommand{\Nd}{\mathbb{N}}

\newcommand{\Rd}{\mathbb{R}}

\newcommand{\eps}{\varepsilon}

\newcommand{\tends}[1]{\xrightarrow[#1]{}}
\newcommand{\tendsd}{\xrightarrow{d}}        % ----> in distribution
\newcommand{\tendsll}{\xrightarrow{L_2}}     % ----> L2 
\newcommand{\tendsas}{\xrightarrow{a.s.}}    % ----> a.s.
\newcommand{\imply}{\Rightarrow}             % implication
\newcommand{\abs}[1]{\lvert #1\rvert}                     % abs value

\let\0=\varnothing
\newcommand{\1}{\mathds{1}}
\newcommand{\2}[1]{\mathds{1}\left\{#1\right\}}

\DeclareMathOperator{\Prob}{P}
\DeclareMathOperator{\Mean}{E}

\DeclareMathOperator{\Var}{Var}

\newcommand{\e}{\mathrm{e}}                  % exponent

\theoremstyle{plain}
\newtheorem{thm}{\bfseries Theorem.}[section]

\newtheorem{prop}{\bfseries Proposition.}[section]
\newtheorem{cor}{\bfseries Corollary.}[section]

\numberwithin{equation}{section}

\renewcommand{\le}{\leqslant}
\renewcommand{\ge}{\geqslant}

\title{On the Optimal Pairwise Group Testing Algorithm}

\author[1]{Ugnė Čižikovienė}
\author[1]{Viktor Skorniakov\thanks{corresponding author; e-mail: viktor.skorniakov@mif.vu.lt}}
\affil[1]{Institute of Applied Mathematics, Faculty of Mathematics and Informatics, Vilnius University, Naugarduko 24, Vilnius LT-03225, Lithuania}

\date{}

\begin{document}

\maketitle

\begin{abstract}
    Originally suggested for the blood testing problem by Dorfman in 1943, an idea of Group Testing (GT) has found many applications in other fields as well. Among many (binomial) GT procedures introduced since then, in 1990, Yao and Hwang proposed the Pairwise Testing Algorithm (PTA) and demonstrated that PTA is the \emph{unique} optimal nested GT procedure provided the probability of contamination lies in $\left[1-\frac{\sqrt{2}}{2},\frac{3-\sqrt{5}}{2}\right]$.
    
    Despite the fundamental nature of the result, PTA did not receive considerable attention in the literature. In particular, even its basic probabilistic properties remained unexplored. In this paper, we fill the gap by providing an exhaustive characterization of probabilistic PTA properties.
\end{abstract}
% derive an explicit representation for the number of tests when applying PTA. The latter representation, given as the sum of weakly dependent indicators, is well suited for further exploration of probabilistic properties of PTA and we illustrate this via examples.

\section{Introduction}\label{s:intro}

% Plan:
% description of the GT idea
% optimal algorithm in terms of the mean number of tests in general 
% class of nested algorithm, Ungar's result and other related results
% no treatment in the literature
% our results briefly and structure of the paper

\emph{Group Testing} (GT) refers to a special kind of technique used to identify defective items in a given set. It is widely applied in very diverse areas. The list includes (but is not limited to) quality control, communication and security networking, genetics, experimental physics, an estimation of parameters from probability models, screening for the infectious diseases like HIV, hepatitis and, very recently, COVID-19 (for domain specific references see e.g. \cite{malinovsky_revisiting_2019}). The main idea underlying the method is as follows. Given a set of items to test, one should replace testing of single items by testing of groups of items. Some of these groups, however, can contain single items as well. For example, to identify defective lights in a given set of strings of lights, one can adopt the following scheme for each string: first test the whole string of lights and then retest each single bulb only in case the whole string does not function properly. This idea was first announced in 1943 by Dorfman \cite{dorfman_detection_1943} who suggested to test pools containing $n$ blood samples and then repeatedly retest only the samples of patients belonging to infected pools. The rationale behind is quite obvious. If the prevalence of the disease is small, then quite often the pooled sample is clean. Hence, instead of testing each sample out of $n$ and consuming this way $n$ test kits, one ends up with a single test. In the literature, the described GT procedure is usually termed as Dorfman scheme, or Dorfman algorithm. Since the appearance of the seminal paper \cite{dorfman_detection_1943}, myriads of other GT algorithms originated. The defining feature of each such algorithm is an average number of tests $E_n$ required to identify all defectives in a set spanning $n$ items. Given $n$, an algorithm achieving minimal possible value of $E_n$ is called an \emph{optimal algorithm}. Characterization of the optimal algorithm without any further assumptions seems to be unfeasible. Therefore, one often operates under the following \emph{Binomial Testing Assumptions} (BTA).
\begin{itemize}
    \item[] \underline{BTA1}: all tested items are independent.
    \item[] \underline{BTA2}: each item is contaminated with the same constant probability $p\in(0,1)$.
    \item[] \underline{BTA3}: pooling does not change operating characteristics (namely, sensitivity and specificity)  of the test kit\footnote{when talking about this assumption, one often says that there is no dilution effect}.
    \item[] \underline{BTA4}: the test kit is perfect, i.e. its sensitivity and specificity are both equal to 1.
\end{itemize}
From now, and till the end of the paper, we assume that BTA hold. In such case, the tested set of items is called the \emph{binomial set} and $E_n$ depends on $p$ as well. Though it might be tempting to conclude that BTA simplify the matter substantially, the forthcoming short account (highlighting fundamental results) demonstrates that the truth is different. 

In 1960, Ungar \cite{Ungar-1960} proved that, for $p>\frac{3-\sqrt{5}}{2}$, irrespectively of value of $n\in\Nd$, an optimal algorithm is to test one-by-one and the minimal value of $E_n=E_n(p)$ is therefore $n$. In 1987, Du and Ko \cite{du_completeness_1987} proved that finding an optimal BTA based algorithm is an NP--complete problem (no polynomial time solution is known \cite{enwiki22-NP-completeness}). In 1988, Yao and Hwang \cite{Yao-monotonicity-1988} proved that $\left(0,\frac{3-\sqrt{5}}{2}\right)\times\{2,3,\ldots\}\ni (p,n)\mapsto E_n(p)$ is monotonically increasing in each argument. Finally, in 1990, the same authors \cite{yao_optimal_1990} proposed the Pairwise Testing Algorithm (PTA) and demonstrated that it is a \emph{unique optimal nested} BTA based algorithm if and only if $p\in\left[1-\frac{\sqrt{2}}{2},\frac{3-\sqrt{5}}{2}\right]$. Nested algorithms are defined by the the following property: if the contaminated subset $C$ is identified, then the next subset to be tested is the proper subset of $C$. Though the optimal nested algorithm is not optimal in the class of all possible GT algorithms, it was demonstrated by Sobel \cite{sobel_group_1960}, \cite{Sobel:1967} that it is nearly optimal over all algorithms. Hence importance of the a fore mentioned result \cite{yao_optimal_1990} on PTA. 

Surprisingly, yet it turns out that an exploration of the properties of PTA did not receive a considerable attention in the literature. Even more, out of 15 citing references \cite{PTA1-abrahams1997code, PTA2-xie2001regression, PTA3-xie2001group,PTA4-tatsuoka2003sequential,PTA5-abrahams1994huffman,PTA6-chi2009optimal,PTA7-golin2001combinatorial,PTA8-li2012nonparametric,PTA9-abrahams1993improved,PTA10-golin2004algorithms,PTA11-malinovsky2020conjectures,PTA12-malinovsky2021nested,PTA13-karimian2005benchmarking,PTA14-ferguson2004optimal,malinovsky_revisiting_2019} retrieved by us\footnote{the list was generated on 28th of June, 2022; non English references were excluded} from \href{https://scholar.google.com/}{Google Scholar}, Malinovsky \cite{PTA11-malinovsky2020conjectures} was the only who investigated a problem having a direct relationship to PTA. All others touched the work of Yao and Hwang \cite{yao_optimal_1990} merely as a reference having a connection to GT with a mild relation to their own problem. These circumstances motivated the current work aiming to give a broader probabilistic characterization of the PTA. Our analysis resulted in the following results for the (properly scaled and/or centered) number of tests performed by PTA: exact analytical expression of the moment generating function (MGF), strong law of large numbers (SLLN), central limit theorem (CLT) and large deviations principle (LDP).

The rest part of the paper is organized as follows. In Section \ref{s:res_and_apps}, we introduce notions and state the previously announced results in detail. Section \ref{s:discussion} contains a short accompanying discussion. Finally, there is an appendix devoted to the proofs.

\section{Results}\label{s:res_and_apps}
We first introduce the PTA by quoting the definition given in \cite{yao_optimal_1990}, Section 2:
\emph{
\begin{itemize}
    \item[(i)] If no contaminated set exists, then always test a pair from the binomial set unless only one item is left, in which case we test that item.
    \item[(ii)] If a contaminated pair is found, test one item of that pair. If that item is good, we deduce that the other is defective. Thus, we classify both items and only a binomial set remains to be classified. If the tested item is defective, the other item together with the remaining binomial set forms a new binomial set.
\end{itemize}}
\noindent Intending to present a full picture, we also restate the main result of Yao and Hwang \cite{yao_optimal_1990}.
\begin{thm}[\cite{yao_optimal_1990}, Theorem 1]
The pairwise testing algorithm is the unique (up to the substitution of equivalent items) optimal nested algorithm for all $n$ if and only if $\frac{2-\sqrt{2}}{2}\le p\le \frac{3-\sqrt{5}}{2}$.
\end{thm}

Let $T_n$ denote the number of conducted tests required for an identification of all defectives in a given binomial set having $n$ items, and let $X_i,i=1,\dots,n$, be an indicator of an $i$th item status. In view of introductory discussion, $X_i\sim Be(p)$ are independent random variables with $p\in(0,1)$ having a meaning of probability of being defective. Also, let $\Bar{X}_i:=1-X_i, q:=1-p$, and
\begin{equation}\label{e:M0_M1}
M_0=\begin{pmatrix}
1 & 0 \\
1 & 0 \\
\end{pmatrix},\quad
M_1=\begin{pmatrix}
0 & 1 \\
1 & 0 \\
\end{pmatrix}.    
\end{equation}
Our first result gives an explicit expression for $T_n$ in terms of the above quantities.
\begin{prop}\label{p:Tn_explicated}
Let $A= \{M_1,M_0M_1\}$ and $B_k = \begin{pmatrix}
X_k & \Bar{X}_k \\
1 & 0 \\
\end{pmatrix}$ for $k=1,\ldots,n$.
Then $T_2=3X_2 + \bar{X}_2(1+X_1), T_3=2 + \bar{X}_3X_2+X_3 T_2$, and
\begin{multline}\label{e:Tn_explicitly}
    T_n = 1+ X_n(\bar{X}_{n-1}X_{n-2}+2) + X_{n-1} +\\ 
    \sum_{j=3}^{n-1}(\bar{X}_{j-1}X_{j-2}+X_{j-1}+1)\left(X_j+\bar{X}_{j-1} \2{B_n B_{n-1}\cdots B_{j+1}\in A}\right) + \\
    X_2+\Bar{X}_2 \2{B_n B_{n-1} \cdots B_{3}\in A}\text{ for }n\ge 4.
\end{multline}
\end{prop}
\noindent The expression above provides insight into the structure of $T_n$ whereas the next one completely characterizes its distribution.

\begin{prop}\label{p:Tn_mgf}
Let $M_{T_n}(\lambda)$ denote the moment generating function of $T_n$ at $\lambda\in\Rd$. Put
\begin{gather}
    \alpha_i=\alpha_i(\lambda)=\frac{1}{2}\left(p\e^{2\lambda} + (-1)^i\sqrt{p^2\e^{4\lambda} +4q\e^\lambda(q+p\e^{\lambda})} \right), \ i=0,1;\label{e:alpha_i}\\
    \kappa_n=\kappa_n(\lambda)=\frac{\alpha_0^n-\alpha_1^n}{\alpha_0-\alpha_1}\text{ for }n\ge 0.\label{e:kappa_n}
\end{gather}
Then 
\begin{multline}\label{e:Tn_mgf}
    M_{T_n}(\lambda)=\e^{2\lambda}\Big[
    \bigl((1-q)^2\e^{3\lambda} + q(1-q)^2\e^{2\lambda} + q(1-q^2)\e^{\lambda} + q^2\bigr)\kappa_{n-2} +\\
    q\left((1-q)^2\e^{3\lambda} +  q(1-q)(2-q)\e^{2\lambda} + 2q^2(1-q)\e^\lambda + q^3 \right)
    \kappa_{n-3}\Big]
\end{multline}
{ for }$n\ge 3$.
\end{prop}
\noindent The remaining results are the consequences of the previous one.
\begin{cor}\label{c:moments}
$\Mean T_n=n\frac{2-q^2}{1+q}+\frac{q^2+q-1}{(1+q)^2}(1-(-q)^n)$,
\begin{multline}\label{e:Var_T_n}
    \Var T_n=\\
    n\frac{\left(1 - q\right)}{\left(q + 1\right)^{3}}\left(q \left(q^{3} + 3 q^{2} + 5 q + 4\right) + \left(- q\right)^{n} \left(2 q + 4\right) \left(q^{2} + q - 1\right)\right) +\\ 
    \frac{\left(1 - \left(- q\right)^{n}\right)}{\left(q + 1\right)^{4}}\left(q \left(5 q^{2} + 3 q - 7\right) + \left(- q\right)^{n} \left(q^{2} + q - 1\right)^{2}\right), n\ge 3.
\end{multline}
\end{cor}

\begin{cor}\label{c:asymptotic_results}
The following asymptotic results apply to $T_n$ as $n\to\infty$.
\begin{itemize}
    \item[] \underline{LLN}: $\frac{T_n}{n}\tendsll \frac{2-q^2}{1+q}$ and $\frac{T_n}{n}\tendsas \frac{2-q^2}{1+q}$.
    \item[] \underline{CLT}: $\sqrt{n}\left(\frac{T_n}{n}-\frac{2-q^2}{1+q}\right)\tendsd N(0,\sigma^2)$, $\sigma^2=\frac{q\left(1 - q\right) \left(q^{3} + 3 q^{2} + 5 q + 4\right)}{\left(q + 1\right)^{3}}$.
    \item[] \underline{LDP}: $\frac{T_n}{n}$ satisfies Large Deviation Principle (LDP) with a good rate function $I$ equal to the Legendre transform of $\Rd\ni\lambda\mapsto \ln\alpha_0(\lambda)$ with $\alpha_0(\lambda)$ given by \eqref{e:alpha_i}. That is, for any closed $C\subset\Rd$ and any open $O\subset\Rd$,
    \begin{gather*}
    \limsup_{n\to\infty}\frac{1}{n}\ln\Prob\left(\frac{T_n}{n}\in C\right)\leq -\inf_{x\in C} I(x)\\
    \text{and}\\
        -\inf_{x\in O} I(x)\leq \liminf_{n\to\infty}\frac{1}{n}\ln\Prob\left(\frac{T_n}{n}\in O\right),
    \end{gather*}
    where $I(x) = \sup_{\lambda\in\Rd}\left(x\lambda - \ln\alpha_0(\lambda)\right)$.
\end{itemize}
\end{cor}

\section{Discussion}\label{s:discussion}

There are several reasons supporting relevance of our analysis.
\begin{itemize}
    \item Though the definition of an optimal algorithm is usually tailored to an average number of tests, when choosing between several  algorithms, it is desirable to evaluate their performance taking into account multiple aspects. For example, an algorithm \emph{A1} may perform slightly better than \emph{A2} in terms of an average number of tests. However, \emph{A1} may have considerably larger variance than \emph{A2} and, therefore, the previously mentioned slight gain of \emph{A1} could be gladly traded by the practitioner in favour of \emph{A2}.
    \item We have already mentioned that the importance of PTA remained unrecognized in the literature and there is more to say on that.
    \begin{itemize}
        \item Many GT algorithms described in the literature (including pioneering Dorfman's algorithm of \cite{dorfman_detection_1943}) have limited applicability due to the \emph{dilution effect}. To be more precise, for a typical algorithm of this kind to perform optimally for a given $p$, one has to test items by grouping them into pools of size $n=n(p)$. If this $n$ is large (say 64 items or even more), the operating characteristics (sensitivity and specificity) of the test kit at hand may become unacceptably low (aka dilute) making this way the algorithm unsuitable for that particular application\footnote{in theory, BTA3 stated in the Introduction prevents from this; however, in practice, it may be a serious obstacle}. With respect to this property, PTA is a very favourable option: it requires only pools of size $n=2$, and this holds true for all $p$'s in the region of its optimality $\left[\frac{2-\sqrt{2}}{2},\frac{3-\sqrt{5}}{2}\right]$.
        \item The region $\left[\frac{2-\sqrt{2}}{2},\frac{3-\sqrt{5}}{2}\right]$ where PTA performs optimally is bounded away from 0 in contrast to many other GT algorithms which do better for $p$'s close to 0. In certain applications this property may be of significant importance. For example, consider a screening for a quite widespread infectious disease.
    \end{itemize}
    \item Yao and Hwang \cite{yao_optimal_1990} conjectured that there exists such $p_0\in\left[\frac{2-\sqrt{2}}{2},\frac{3-\sqrt{5}}{2}\right]$ that for $p\in\left[p_0,\frac{3-\sqrt{5}}{2}\right]$ PTA is optimal over all (not necessarily nested) algorithms satisfying BTA.
    \item Our Prop. \ref{p:Tn_explicated} demonstrates that, despite apparently simple recurrence governing evolution of $T_n$ (see Eq. \eqref{e:Tn_reccurence}), the resulting dependence structure is not so simple. At least we were not able to analyze its behaviour neither by making use of Markov chains theory, nor by making use of martingale theory. A well developed apparatus of weakly dependent sequences also did not promise easy deduction of Corollary \ref{c:asymptotic_results}. More than that, even direct moment calculation exercise, though accomplishable for $\Mean T_n$ at a reasonable price (see Lemma in Section 4 of \cite{yao_optimal_1990}), becomes much more involved when it comes to $\Var T_n$ and higher order moments. This way, $(T_n)_{n\ge 2}$ yields an example of a sequence of positive integer valued random variables having an interesting probabilistic structure encountered in practical application and not designed artificially for learning or other purposes.
\end{itemize}
In view  of the said above, our input seems to be plausible. Moreover, we are inclined to think that it may be useful for the solution of a couple of unresolved conjectures. Namely, the one stated by Yao and Hwang in \cite{yao_optimal_1990} and mentioned above, and the generalized PTA optimality conjecture stated in \cite{PTA11-malinovsky2020conjectures}.

\appendix
\section{Proofs}
\emph{Proof of Proposition \ref{p:Tn_explicated}}. By the description of the testing procedure,
\begin{multline}\label{e:Tn_reccurence}
    T_n = (1+T_{n-2})\1\{X_n + X_{n-1}=0 \}+(2+T_{n-1})\1\{X_n + X_{n-1}>0\}X_n + \\
    (2+T_{n-2})\1\{X_n+X_{n-1}>0\}\Bar{X}_n=\Big[
    \text{since } \1\{X_n + X_{n-1}=0\} = \Bar{X}_n\Bar{X}_{n-1}, \\
    \1\{X_n + X_{n-1}>0\}X_n= X_n,\quad \1\{X_n + X_{n-1}>0\}\Bar{X}_n= \Bar{X}_nX_{n-1}
    \Big]=\\
    \Bar{X}_n(1+X_{n-1}) +2X_n + T_{n-1}X_n + T_{n-2}\Bar{X}_n.
\end{multline}
Put
\begin{equation*}
    t_1 = \begin{pmatrix}
    1\\
    0
    \end{pmatrix},\quad
    A_1 = t_1,\quad
    t_k \stackrel{k\ge 2}{=}
    \begin{pmatrix}
    T_k\\
    T_{k-1}
    \end{pmatrix},\quad
    A_k \stackrel{k\ge 2}{=}
    \begin{pmatrix}
    \bar{X}_k(1+X_{k-1}) + 2X_k\\
    0
    \end{pmatrix},
\end{equation*}
and let $B_k,k\ge 1$, be as in the statement of the Proposition. From \eqref{e:Tn_reccurence} it follows that
\begin{multline*}
    t_n = A_n + B_n t_{n-1} = \ldots = A_n + \sum_{k=1}^{n-2} B_n B_{n-1}\dots B_{n-k+1}A_{n-k} + B_n\dots B_{2}t_1 =\\
    A_n + \sum_{j=3}^n B_n\dots B_jA_{j-1} + B_n\dots B_{2}t_1=A_n + \sum_{j=2}^n B_n\dots B_jA_{j-1}.
\end{multline*}
Let $M_0,M_1$ be given by \eqref{e:M0_M1}. Denoting
\begin{equation}\label{e:M2_M3}
    M_2 = \begin{pmatrix}
    1 & 0\\
    0 & 1\\
    \end{pmatrix},\quad
    M_3 = \begin{pmatrix}
    0 & 1\\
    0 & 1\\
    \end{pmatrix},\quad S=\{M_0,M_1,M_2,M_3\},
\end{equation}
we have that $S$ forms a semi-group with respect to ordinary matrix multiplication since
\begin{gather}
    M_0^2 = M_0,\quad M_0M_1 = M_3,\quad M_0M_3 = M_3,\quad M_1M_0 = M_0,\quad M_1^2 = M_2,\nonumber\\
    M_1M_3 = M_3,\quad M_3M_0 = M_0,\quad M_3M_1 = M_0,\quad M_3^2 = M_3.\label{e:M_closure_props}
\end{gather}
Let $J_i = \{j\in\{2,,\dots,n\}\mid X_j=i\},i=0,1$. Note that $\forall\, i\ B_i=X_iM_0 + \bar{X}_iM_1\in S$ and that $M_0$ is an absorbing element of $S$. Therefore, by \eqref{e:M_closure_props}
\begin{equation*}
    \sum_{j\in J_1} B_n\dots B_j A_{j-1}= \sum_{j\in J_1} M_0A_{j-1}=
    \sum_{j=2}^n X_jM_0 A_{j-1}
\end{equation*}
and
\begin{multline*}
    \sum_{j\in J_0} B_n\dots B_j A_{j-1}=
    \sum_{j=2}^n \bar{X}_j\Big(
        \2{B_n\dots B_j M_1=M_0}M_0 +\\ \2{B_n\dots B_j M_1=M_1}M_1 +
        \2{B_n\dots B_j M_1=M_2}M_2 +\\ \2{B_n\dots B_j M_1=M_3}M_3
    \Big) A_{j-1}.
\end{multline*}
To extract $T_n$ from $t_n$, it suffices to multiply $t_n$ by $(1\ 0)$ from the left. Since 
\begin{equation*}
    (1\ 0)M_iA_{j-1}=\begin{cases}
        &0, \text{ for }i = 1,3\text{ and all }j;\\
        &\bar{X}_{j-1}(1+X_{j-2}) + 2X_{j-1}, \text{ for }i=0,2\text{ and }j\ge 3;\\
        &1,\text{ for }i=0,2\text{ and }j=2,
    \end{cases}
\end{equation*}
after the collection of terms, we finally end up with an expression \eqref{e:Tn_explicitly}. $\qed$

\smallskip\noindent\emph{Proof of Proposition \ref{p:Tn_mgf}. \underline{Step 1: auxiliary recurrence.}} For $\lambda_1,\lambda_2\in\Rd$, let 
\begin{equation}\label{e:cond_expectations}
    M_{i,n}(\lambda_1,\lambda_2)=\Mean\left(\e^{\lambda_1 T_{n} + \lambda_2 T_{n-1}} \mid X_n=i\right),\quad i=0,1.
\end{equation}
By equation \eqref{e:Tn_reccurence},
\begin{align}\label{e:Mn_reccurence}
    M_{0,n}(\lambda_1,\lambda_2) & =\Mean\left(\e^{\lambda_1(1+X_{n-1} + T_{n-2}) + \lambda_2 T_{n-1}}\right)= \nonumber\\
    &p\Mean\left(\e^{\lambda_1(2+ T_{n-2}) + \lambda_2 T_{n-1}}\mid X_{n-1}=1\right) + \nonumber\\
    &q\Mean\left(\e^{\lambda_1(1+ T_{n-2}) + \lambda_2 T_{n-1}}\mid X_{n-1}=0\right)=\nonumber\\
    &p\e^{2\lambda_1} M_{1,n-1}(\lambda_2,\lambda_1) + 
    q\e^{\lambda_1} M_{0,n-1}(\lambda_2,\lambda_1);\nonumber\\
    M_{1,n}(\lambda_1,\lambda_2) &= \Mean\left(\e^{\lambda_1(2+ T_{n-1}) + \lambda_2 T_{n-1}}\right)=\nonumber\\
    &\e^{2\lambda_1}\left(p M_{1,n-1}(\lambda_1+\lambda_2,0) + q M_{0,n-1}(\lambda_1+\lambda_2,0)\right).
\end{align}
For $\lambda\in\Rd$, put 
\begin{multline}\label{e:matrices}
    m_{1,n}=m_{1,n}(\lambda)=M_{1,n}(\lambda,0),\quad m_{2,n}=m_{2,n}(\lambda)=M_{0,n}(\lambda,0),\\
    m_{3,n}=m_{3,n}(\lambda)=M_{1,n}(0,\lambda),\quad m_{4,n}=m_{4,n}(\lambda)=M_{0,n}(0,\lambda);\\
    A = A(\lambda) = \e^{2\lambda}\begin{pmatrix}
    p & q \\
    0 & 0\\
    \end{pmatrix},\quad
    B = B(\lambda) = \begin{pmatrix}
    0 & 0\\
    \e^{2\lambda}p & \e^{\lambda}q \\
    \end{pmatrix},\\
    C = \begin{pmatrix}
    p & q\\
    p & q\\
    \end{pmatrix},\quad
    O = \begin{pmatrix}
    0 & 0\\
    0 & 0\\
    \end{pmatrix}.
\end{multline}
From \eqref{e:Mn_reccurence} it then follows that $m_n=m_n(\lambda)=(m_{1,n},m_{2,n},m_{3,n},m_{4,n})^\top$ satisfies recurrent equation
\begin{equation}\label{e:mnRecurrence}
    m_n = 
    \begin{pmatrix}
    A & B\\
    C & O\\
    \end{pmatrix}m_{n-1}=\dots = 
    \begin{pmatrix}
    A & B\\
    C & O\\
    \end{pmatrix}^{n-1}m_1
\end{equation}
Writing
\[
\begin{pmatrix}
    A & B\\
    C & O\\
    \end{pmatrix}^{n}=\begin{pmatrix}
    A_n & B_n\\
    C_n & D_n\\
    \end{pmatrix}
\] and applying inductive argument, one finds out that the $2\times2$ blocks $A_n,B_n,C_n,D_n$ satisfy 
\begin{align}
    &\begin{cases}
    &A_n = A A_{n-1} + B C_{n-1},\\
    &C_n = C A_{n-1};\\
    \end{cases} \label{e:block_system1} \\
    &\begin{cases}
    &B_n = A B_{n-1} + B D_{n-1},\\
    &D_n = C B_{n-1};\\
    \end{cases} \label{e:block_system2}
\end{align}
with $A_0 = D_0 = Id$ and $ C_0=B_0=O$. Consider system \eqref{e:block_system1}. Since $A=\begin{pmatrix}\e^{2\lambda} & 0\\0 & 0\\\end{pmatrix}C$, we have that
\begin{multline}\label{e:An_recurrence}
    A_n = \begin{pmatrix}\e^{2\lambda} & 0\\0 & 0\\\end{pmatrix} CA_{n-1} + BC_{n-1}=
    \begin{pmatrix}\e^{2\lambda} & 0\\0 & 0\\\end{pmatrix}C_{n} + B C_{n-1}.
\end{multline}
Therefore,
\begin{equation}\label{e:Cn_recurrence}
    C_n = C\left(\begin{pmatrix}\e^{2\lambda} & 0\\0 & 0\\\end{pmatrix}C_{n-1} + B C_{n-2}\right).
\end{equation}
Let $\kappa_n$ be defined by \eqref{e:kappa_n}. We claim that $C_n=\kappa_n C$ solves \eqref{e:Cn_recurrence}. For $n=2$ (as well as $n=0,1$) the claim holds by direct check. Assume it holds for $k\le n$ with $n\ge 2$. Applying inductive assumption and multiplying,
\begin{multline*}
    C_{n+1}=C\left(\kappa_n\begin{pmatrix}\e^{2\lambda} & 0\\0 & 0\\\end{pmatrix} C + \kappa_{n-1} BC\right)=\\
    \left[ 
    C\begin{pmatrix}\e^{2\lambda} & 0\\0 & 0\\\end{pmatrix}C = p\e^{2\lambda}C,\
    C B C=q\e^\lambda(q+p\e^\lambda)C
    \right]=\\
    (p\e^{2\lambda} \kappa_{n}+q\e^\lambda(q+p\e^\lambda)\kappa_{n-1})C=\kappa_{n+1}C
\end{multline*}
since an expression for $\kappa_n$ given in \eqref{e:kappa_n} is precisely the solution of the second order linear difference equation
\begin{equation*}
    \kappa_{n+1}=p\e^{2\lambda} \kappa_{n}+q\e^\lambda(q+p\e^\lambda)\kappa_{n-1},\quad\kappa_1=1,\quad\kappa_0=0.
\end{equation*}
Substituting $C_n=\kappa_n C$ to \eqref{e:An_recurrence}, we obtain an expression for $A_n$.

System \eqref{e:block_system2} is handled in the same way by noting that it is identical to \eqref{e:block_system1} and only the initial conditions differ leading thereby to the following solution: 
\begin{equation}\label{e:B_n_Dn_solution}
    D_n= \kappa_{n-1}D_2,\quad B_n = \begin{pmatrix}\e^{2\lambda} & 0\\0 & 0\\\end{pmatrix}D_{n} + B D_{n-1}\text{ for }n\ge 1.
\end{equation}

\underline{\emph{Step 2: final expression.}} From the results of \emph{Step 1}, we obtain an expression for $m_n$ given by \eqref{e:mnRecurrence} since $m_1$ is readily available and equal to\footnote{note that $T_1\equiv 1, T_0\equiv0$} $(\e^{\lambda},\e^\lambda,1,1)^\top$:
\begin{equation*}
    m_n = \begin{pmatrix}
    (\e^\lambda A_{n-1} + B_{n-1})\begin{pmatrix}
    1\\
    1\\
    \end{pmatrix}\\
    (\e^\lambda C_{n-1} + D_{n-1})\begin{pmatrix}
    1\\
    1\\
    \end{pmatrix}\\
    \end{pmatrix}.
\end{equation*}
Noting that 
\begin{equation*}
    \Mean \e^{\lambda T_n}=p\Mean \left(\e^{\lambda T_n}\mid X_n=1\right) + 
    q\Mean \left(\e^{\lambda T_n}\mid X_n=0\right)=p m_{1,n} + q m_{2,n},
\end{equation*}
we finally arrive to expression \eqref{e:Tn_mgf}. $\qed$

\smallskip\noindent\emph{Proof of Corollary \ref{c:moments}}. Recall that the $k$-th derivative of the moment generating function evaluated at 0 yields the $k$-th moment. Therefore, to obtain the announced formulae, one simply needs to differentiate expression \eqref{e:Tn_mgf}. Though conceptually an exercise is trivial, the calculations require tedious work. Therefore, we provide key steps and some intermediate quantities yet omit the detailed listing in order not to overwhelm the paper with the trivial content. For the sake of convenience, we make change of variables $x=\e^{\lambda}$ and work with probability generating function $G(x)=\Mean x^{T_n}=M_{T_n}(\ln\lambda)$. By \eqref{e:kappa_n}--\eqref{e:Tn_mgf} and slight abuse of notation,
\begin{align}
    &G(x)= g_1(x)\kappa_{n-2}(x) + g_2(x)\kappa_{n-3}(x)\text{ with} \nonumber\\
    &g_1(x)=\bigl((1-q)^2x^{3} + q(1-q)^2x^{2} + q(1-q^2)x + q^2\bigr)x^2,\nonumber\\
    &g_2(x)=q\left((1-q)^2x^{3} +  q(1-q)(2-q)x^{2} + 2q^2(1-q)x + q^3 \right)x^2,\nonumber\\
    &\alpha_i=\alpha_i(x)=\frac{1}{2}\left(px^{2} + (-1)^i\sqrt{p^2x^{4} +4qx(q+px)} \right), \text{ for } i=0,1,\text{ and }\nonumber\\
    &\kappa_n=\kappa_n(x)=\frac{\alpha_0^n(x)-\alpha_1^n(x)}{\alpha_0(x)-\alpha_1(x)}\text{ for }n\ge 0.\label{e:pgf1}
\end{align}
Then 
\begin{multline}\label{e:ETn_via_derivative}
    \Mean T_n =G^\prime(1)=
    g_1^\prime(1)\kappa_{n-2}(1) + g_1(1)\kappa_{n-2}^\prime(1) + \\
    g_2^\prime(1)\kappa_{n-3}(1) + g_2(1)\kappa_{n-3}^\prime(1)
\end{multline}
and
\begin{multline}\label{e:ETnTn-1_via_derivative}
    \Mean T_n(T_n-1)=G^{\prime\prime}(1)=
    g_1^{\prime\prime}(1)\kappa_{n-2}(1) + 2g_1^\prime(1)\kappa_{n-2}^\prime(1) + g_1(1)\kappa_{n-2}^{\prime\prime}(1)+\\
    g_2^{\prime\prime}(1)\kappa_{n-3}(1)+ 2g_2^\prime(1)\kappa_{n-3}^\prime(1)+g_2(1)\kappa_{n-3}^{\prime\prime}(1).
\end{multline}
Therefore, $\Var T_n=G^{\prime\prime}(1)+G^{\prime}(1)-\left(G^{\prime}(1)\right)^2$ and to verify the announeced formulae, one needs to check the validity of the equalities
\begin{align*}
    &\alpha_0(1)=1,\quad \alpha_1(1)=-q,\quad\alpha_0^\prime(1)=\frac{2-q^2}{1+q},\quad \alpha_1^\prime(1)=-\frac{q^2}{1+q},\\
    &\alpha_0^{\prime\prime}(1)=4\frac{1 - q}{q + 1} - \frac{2}{\left(q + 1\right)^{3}},\quad
    \alpha_1^{\prime\prime}(1)= -\frac{2 \left(1 - q\right)^{2}}{q + 1} + \frac{2}{\left(q + 1\right)^{3}},\\
\end{align*}
\begin{align*}
    &g_1(1)=1,\quad g_2(1)=q,\\ 
    &g_1^\prime(1)=q^3-q^2-3q+5,\quad g_2^\prime(1)=-q\left(q^2+2q-5 \right),\\
    &g_1^{\prime\prime}(1)=6 q^{3} - 2 q^{2} - 22 q + 20,\quad g_2^{\prime\prime}(1)=2 q \left(q^{3} - 2 q^{2} - 8 q + 10\right),\\
\end{align*}
\begin{align*}
    &\kappa_n(1)=\frac{1-(-q)^n}{1+q},\quad \kappa_n^\prime(1)=n\frac{2-q^2}{(1+q)^2} + \frac{(-q)^n(2-q(1+q)n)-2}{(1+q)^3},\\
    &\kappa_n^{\prime\prime}(1)=\frac{2 n \left(1 - q\right) \left(2+\left(1 - q\right)\left(- q\right)^{n - 1} \right)}{\left(q + 1\right)^{2}} +\\ 
    &\frac{n \left(n - 1\right) \left(\left(2 - q^{2}\right)^{2} - \left(- q\right)^{n + 2} \right) - 
    2\left(1 - q\right) \left(3 - q\right)\left(1 - \left(- q\right)^{n}\right)}{\left(q + 1\right)^{3}}-\\
    &\frac{2 n \left(- 2 q^{2} + 5 + \left(- q\right)^{n - 1} \left(2 q^{2} + 1\right) \right)}{\left(q + 1\right)^{4}} + 12\frac{1 - \left(- q\right)^{n}}{\left(q + 1\right)^{5}},
\end{align*}
plug them into \eqref{e:ETn_via_derivative}--\eqref{e:ETnTn-1_via_derivative}, and carefully collect the terms. $\qed$

\smallskip\noindent\emph{Proof of corollary \ref{c:asymptotic_results}. \underline{Step 1: expansions}}. Applying Taylor's formula, we obtain the following equalities (for $\lambda\to 0$):
\begin{align}\label{e:primary_expansions}
    &p^2\e^{4\lambda} + 4q\e^\lambda\left(q+p\e^\lambda\right) = (1+q)^2\left[ 
    1 + \frac{4\lambda}{(1+q)^2} + 2\lambda^2\left(\frac{2-q}{1+q}\right)^2 + O(\lambda^3)
    \right];\nonumber\\
    &\sqrt{
        p^2\e^{4\lambda}+4q\e^{\lambda}(q+p\e^\lambda)
    } = 1 + \frac{2\lambda}{(1+q)^2} +\nonumber\\
    &\lambda^2\left(\left(\frac{2-q}{1+q}\right)^2 - \frac{2}{(1+q)^4}\right) + O(\lambda^3);\nonumber\\
    &\alpha_0 = 1 + \lambda\frac{2-q^2}{1+q} + \frac{\lambda^2}{2}\left(
    2(1-q) + \frac{(2-q)^2}{(1+q)} - \frac{2}{(1+q)^3}
    \right) + O(\lambda^3); \nonumber\\
    &\alpha_1 = -q - \lambda\frac{q^2}{1+q} + \frac{\lambda^2}{2}\left(
    2(1-q) - \frac{(2-q)^2}{(1+q)} + \frac{2}{(1+q)^3}
    \right) + O(\lambda^3).
\end{align}
Let $c_{ij}$ denote a coefficient near $\lambda^j$ in the expansion of $\frac{\alpha_i}{(-q)^i}$ for $j=0,1,2$ and $i=0,1$. Then
\begin{equation}\label{e:alpha_i_expansions}
    \ln\left(\frac{\alpha_i}{(-q)^i}\right)^n = n\left(
        c_{i1}\lambda + \left( c_{i2} - \frac{c_{i1}^2}{2} \right)\lambda^2
    \right) + O(n\lambda^3).
\end{equation}
Consequently, 
\begin{multline}\label{e:kappa_n_expansion}
    (\alpha_0-\alpha_1)\kappa_n(\lambda) = \alpha_0^n - \alpha_1^n = \e^{n\ln\alpha_0} - (-1)^n\e^{n\ln \left(q\frac{\alpha_1}{-q}\right)} = \\
    \exp\left\{
        n\left(
        c_{01}\lambda + \left( c_{02} - \frac{c_{01}^2}{2} \right)\lambda^2
    \right) + O(n\lambda^3)
    \right\}- \\
    (-1)^n
    \exp\left\{
        n\left(
        c_{11}\lambda + \left( c_{12} - \frac{c_{11}^2}{2} \right)\lambda^2
    \right) + n\ln q + O(n\lambda^3)
    \right\}.
\end{multline}
Finally, let $g_i(x)$ denote the same polynomials as given in \eqref{e:pgf1}. Taylor expanding yields
\begin{equation*}
    g_1(\e^\lambda) = 1 + O(\lambda),\quad g_2(\e^\lambda) = q + O(\lambda).
\end{equation*}
Combining all above, we then obtain the following asymptotic expansion for the moment generating function:
\begin{equation}
    M_{T_n}(\lambda) = \frac{1}{1+q+O(\lambda)}\left(
    (1 + O(\lambda))\kappa_{n-2}(\lambda) + (q + O(\lambda))\kappa_{n-3}(\lambda)
    \right),\quad \lambda\to 0,
\end{equation}
with asymptotic expressions for $\kappa_{n-2},\kappa_{n-3}$ stemming from \eqref{e:kappa_n_expansion}.

\emph{\underline{Step 2: LLN}}. To prove relationship $\frac{T_n}{n}\tendsll \frac{2-q^2}{1+q}$, note that, by Corollary \ref{c:moments},
\begin{multline*}
    \Mean\left(
    \frac{T_n}{n} - \frac{2-q^2}{1+q}
    \right)^2 = \Mean\left(
    \frac{T_n}{n} - \Mean \frac{T_n}{n}
    \right)^2 + \Mean\left(
    \Mean \frac{T_n}{n} - \frac{2-q^2}{1+q}
    \right)^2 =\\
    \frac{1}{n^2}\left(
        \Var T_n + \left( \frac{q^2+q-1}{(1+q)^2}\left(1-(-q)^n\right)\right)^2
    \right) = O\left(\frac{1}{n}\right).
\end{multline*}
To prove a.s. convergence, we show that the following sufficient condition holds:
\begin{equation}\label{e:as_conv_cond}
    \forall \eps>0\quad \sum_{n=2}^\infty \Prob\left( \Big|\frac{T_n}{n} - \frac{2-q^2}{1+q}\Big|>\eps \right)<\infty.
\end{equation}
To this end, we bound the probability
\begin{multline*}
    \Prob\left( \Big|\frac{T_n}{n} - \frac{2-q^2}{1+q}\Big|>\gamma\frac{\ln n}{\sqrt{n}} \right) = \\
    \Prob\left(\frac{T_n}{\sqrt{n}} - \sqrt{n}\frac{2-q^2}{1+q}>\gamma\ln n\right) + 
    \Prob\left(\frac{T_n}{\sqrt{n}} - \sqrt{n}\frac{2-q^2}{1+q}<-\gamma\ln n\right), 
\end{multline*}
where $\gamma>0$ is arbitrary yet fixed constant. By Markov's inequality,
\begin{multline*}
    \Prob\left(\frac{T_n}{\sqrt{n}} - \sqrt{n}\frac{2-q^2}{1+q}>\gamma\ln n\right)\leq\\
    \e^{-\sqrt{n}\frac{2-q^2}{1+q}-\gamma\ln n}\Mean{\e^{\frac{T_n}{\sqrt{n}}}}=
    \e^{-\sqrt{n}\frac{2-q^2}{1+q}-\gamma\ln n}M_{T_n}\left(\frac{1}{\sqrt{n}}\right).
\end{multline*}
From results obtained in Step 1 and after some rearrangement, it follows that
\begin{multline*}
    M_{T_n}\left(\frac{1}{\sqrt{n}}\right)=
    \left(1+O\left(\frac{1}{\sqrt{n}}\right)\right)\Big(
    \e^{c_{01}\sqrt{n} + c_{02} - \frac{c_{01}^2}{2} + O(\frac{1}{\sqrt{n}})} -\\ 
    (-1)^n
    \e^{c_{11}\sqrt{n} + c_{12} - \frac{c_{11}^2}{2} + n\ln q + O(\frac{1}{\sqrt{n}})}\Big).
\end{multline*}
Since, 
\begin{multline*}
    c_{01} = \frac{2-q^2}{1+q}\text{ and } c_{11}\sqrt{n} - \frac{2-q^2}{1+q}\sqrt{n} + c_{12}-\frac{c_{11}^2}{2} +  n\ln q =\\ 
    n\ln q\left( 1+ O\left(\frac{1}{\sqrt{n}}\right)\right),
\end{multline*}
we obtain that
\begin{equation*}
    \e^{-\sqrt{n}\frac{2-q^2}{1+q}-\gamma\ln n}M_{T_n}\left(\frac{1}{\sqrt{n}}\right) = \e^{-\gamma\ln n}O(1) \leq \frac{C_q}{n^\gamma}
\end{equation*}
for some constant $C_q\in(0,\infty)$ independent of $\gamma$. In the same way,
\begin{multline*}
    \Prob\left( \Big|\frac{T_n}{n} - \frac{2-q^2}{1+q}\Big|<-\gamma\frac{\ln n}{\sqrt{n}} \right) \leq \e^{\sqrt{n}\frac{2-q^2}{1+q}-\gamma\ln n} M_{T_n}\left(-\frac{1}{\sqrt{n}}\right)\leq \frac{C_q}{n^\gamma},
\end{multline*}
provided $C_q$ in the previous inequality was chosen large enough. Taking $\gamma>1$, we then have that
\begin{equation*}
    \sum_{n=2}^\infty \Prob\left( \Big|\frac{T_n}{n} - \frac{2-q^2}{1+q}\Big|>\gamma\frac{\ln n}{\sqrt{n}} \right)\leq 2C_q\sum_{n=1}^\infty\frac{1}{n^\gamma} <\infty.
\end{equation*}
Hence \eqref{e:as_conv_cond} and the claim.

\emph{\underline{Step 3: CLT}}. It suffices to show that
\begin{equation*}
    M_{\sqrt{n}\left( \frac{T_n}{n} - \frac{2-q^2}{1+q}\right)}(t)\tends{n\to\infty}
    M_\xi(t),\quad \xi\sim N(0,\sigma^2)
\end{equation*}
for some fixed $\eps>0$ and any fixed $t\in(-\eps,\eps)$. Applying expansions obtained in the \emph{Step 1} and the reasoning similar to that of \emph{Step 2}, we have that
\begin{multline*}
    M_{\sqrt{n}\left( \frac{T_n}{n} - \frac{2-q^2}{1+q}\right)}(t) = \e^{-t\sqrt{n}\frac{2-q^2}{1+q}}
    M_{T_n}\left(\frac{t}{\sqrt{n}}\right) =\\
    \left(1+O\left(\frac{1}{\sqrt{n}}\right)\right)\Big(
    \e^{tc_{01}\sqrt{n} + t^2\left(c_{02} - \frac{c_{01}^2}{2}\right) + O(\frac{1}{\sqrt{n}})} -\\ 
    (-1)^n
    \e^{tc_{11}\sqrt{n} + t^2\left(c_{12} - \frac{c_{11}^2}{2}\right) + n\ln q + O(\frac{1}{\sqrt{n}})}\Big)=\\
    \left(1+O\left(\frac{1}{\sqrt{n}}\right)\right)
    \e^{t^2\left(c_{02} - \frac{c_{01}^2}{2}\right) + O(\frac{1}{\sqrt{n}})} + 
    O\left(q^{n}\right)\tends{n\to\infty}\e^{t^2\left(c_{02} - \frac{c_{01}^2}{2}\right)}.
\end{multline*}
Direct calculations show that $c_{01}-\frac{c_{02}^2}{2}=\frac{\sigma^2}{2}$.

\emph{\underline{Step 4: LDP}}. To prove the final claim, we apply Gärtner-Ellis (GE) Theorem (see \cite{dembo2009large}, Section 2.3) to $Z_n=\frac{T_n}{n}$. First, note that, for any fixed $\lambda\in\Rd$, 
\begin{multline*}
    \alpha_0(\lambda)>\abs{\alpha_1(\lambda)}\imply \lim_{n\to\infty}\frac{\kappa_{n-3}(\lambda)}{\kappa_{n-2}(\lambda)} = \frac{1}{\alpha_0(\lambda)} \imply\\
    \Lambda(\lambda) := \lim_{n\to\infty} \frac{1}{n}\ln M_{Z_n}(n\lambda) = \lim_{n\to\infty} \frac{1}{n}\ln M_{T_n}(\lambda) = \\
    \lim_{n\to\infty} \frac{1}{n}\ln \alpha_{0}^n(\lambda) = \ln\alpha_0(\lambda)\in\Rd.
\end{multline*}
Since $\Rd\ni \lambda\mapsto \Lambda(\lambda)$ is differentiable at every $\lambda\in\Rd$, it follows that all GE assumptions hold and $T_n$ satisfies LDP with a good rate function $I$ equal to the Legendre transform of $\Lambda$. $\qed$

% \begin{multline*}
%     \Var T_n=\\
%     n\frac{\left(1 - q\right)}{\left(q + 1\right)^{3}}\left(q \left(q^{3} + 3 q^{2} + 5 q + 4\right) + \left(- q\right)^{n} \left(2 q + 4\right) \left(q^{2} + q - 1\right)\right) +\\ 
%     \frac{\left(1 - \left(- q\right)^{n}\right)}{\left(q + 1\right)^{4}}\left(q \left(5 q^{2} + 3 q - 7\right) + \left(- q\right)^{n} \left(q^{2} + q - 1\right)^{2}\right)
% \end{multline*}
% \begin{multline}
%     \Mean T_n=\frac{\d}{\d x}G(x)\Big|_{x=1}= 
%     \text{ and } \Mean T_n(T_n-1)=\frac{\d^2}{\d x^2}G(x)\Big|_{x=1}=,
% \end{multline}
% \begin{align*}
% B_k = \begin{pmatrix}
% X_k & \Bar{X}_k \\
% 1 & 0 \\
% \end{pmatrix}\text{ for }k=1,\ldots,n.    
% \end{align*}
% , A_1=(1\ 0)^T,
% &A_k = \begin{pmatrix}\Bar{X}_k(1+X_{k-1})+2X_k \\ 0 \end{pmatrix},\text{ for }k=2,\ldots,n,\text{ and }\\

\bibliographystyle{plain}

\begin{thebibliography}{10}

\bibitem{PTA9-abrahams1993improved}
Julia Abrahams.
\newblock An improved lower bound on the minimum expected number of binomial
  group tests.
\newblock {\em Probability in the Engineering and Informational Sciences},
  7(1):121--124, 1993.

\bibitem{PTA5-abrahams1994huffman}
Julia Abrahams.
\newblock Huffman-type codes for infinite source distributions.
\newblock In {\em Proceedings of IEEE Data Compression Conference (DCC'94)},
  pages 83--89. IEEE, 1994.

\bibitem{PTA1-abrahams1997code}
Julia Abrahams.
\newblock Code and parse trees for lossless source encoding.
\newblock {\em Proceedings. Compression and Complexity of SEQUENCES 1997 (Cat.
  No. 97TB100171)}, pages 145--171, 1997.

\bibitem{PTA6-chi2009optimal}
Xiao-Fei Chi, Xiang-Yang Lou, Mark~CK Yang, and Qing-Yao Shu.
\newblock An optimal dna pooling strategy for progressive fine mapping.
\newblock {\em Genetica}, 135(3):267--281, 2009.

\bibitem{dembo2009large}
A.~Dembo and O.~Zeitouni.
\newblock {\em Large Deviations Techniques and Applications}.
\newblock Stochastic Modelling and Applied Probability. Springer Berlin
  Heidelberg, 2009.

\bibitem{dorfman_detection_1943}
R.~Dorfman.
\newblock The detection of defective members of large populations.
\newblock {\em The Annals of Mathematical Statistics}, 14(4):436--440, 1943.

\bibitem{du_completeness_1987}
Ding-Zhu Du and Ker-I Ko.
\newblock Some {Completeness} {Results} on {Decision} {Trees} and {Group}
  {Testing}.
\newblock {\em SIAM Journal on Algebraic Discrete Methods}, 8(4):762--777,
  October 1987.

\bibitem{PTA14-ferguson2004optimal}
Thomas~S Ferguson and Curtis Tatsuoka.
\newblock An optimal strategy for sequential classification on partially
  ordered sets.
\newblock {\em Statistics \& probability letters}, 68(2):161--168, 2004.

\bibitem{PTA7-golin2001combinatorial}
Mordecai~J Golin.
\newblock A combinatorial approach to golomb forests.
\newblock {\em Theoretical Computer Science}, 263(1-2):283--304, 2001.

\bibitem{PTA10-golin2004algorithms}
Mordecai~J Golin and Kin~Keung Ma.
\newblock Algorithms for infinite huffman-codes.
\newblock In {\em Proceedings of the fifteenth annual ACM-SIAM symposium on
  Discrete algorithms}, pages 758--767, 2004.

\bibitem{PTA13-karimian2005benchmarking}
SA~Mohsen Karimian and Anthony~G Straatman.
\newblock Benchmarking of a 3d, unstructured, finite volume code for
  incompressible navier-stokes equation on a cluster of distributed-memory
  computers.
\newblock In {\em 19th International Symposium on High Performance Computing
  Systems and Applications (HPCS'05)}, pages 11--16. IEEE, 2005.

\bibitem{PTA8-li2012nonparametric}
Mingyu Li and Minge Xie.
\newblock Nonparametric and semiparametric regression analysis of group testing
  samples.
\newblock {\em International Journal of Statistics in Medical Research},
  1(1):60--72, 2012.

\bibitem{PTA11-malinovsky2020conjectures}
Yaakov Malinovsky.
\newblock Conjectures on optimal nested generalized group testing algorithm.
\newblock {\em Applied Stochastic Models in Business and Industry},
  36(6):1029--1036, 2020.

\bibitem{malinovsky_revisiting_2019}
Yaakov Malinovsky and Paul~S. Albert.
\newblock Revisiting {Nested} {Group} {Testing} {Procedures}: {New} {Results},
  {Comparisons}, and {Robustness}.
\newblock {\em The American Statistician}, 73(2):117--125, April 2019.

\bibitem{PTA12-malinovsky2021nested}
Yaakov Malinovsky and Paul~S Albert.
\newblock Nested group testing procedures for screening.
\newblock {\em arXiv preprint arXiv:2102.03652}, 2021.

\bibitem{sobel_group_1960}
M.~Sobel.
\newblock Group {Testing} to {Classify} {Efficiently} all {Defectives} in a
  {Binomial} {Sample}.
\newblock In R.~E. Machol, editor, {\em Information and {Decision}
  {Processes}}, pages 127--161. New York, McGraw Hill, 1960.

\bibitem{Sobel:1967}
M.~Sobel.
\newblock Optimal group testing.
\newblock In {\em Proceedings of the Colloquium on Information Theory}, pages
  411--488, Debrecen (Hungary), 1967. Organized by the Bolyai Mathematical
  Society.

\bibitem{PTA4-tatsuoka2003sequential}
Curtis Tatsuoka and Thomas Ferguson.
\newblock Sequential classification on partially ordered sets.
\newblock {\em Journal of the Royal Statistical Society: Series B (Statistical
  Methodology)}, 65(1):143--157, 2003.

\bibitem{Ungar-1960}
Peter Ungar.
\newblock The cutoff point for group testing.
\newblock {\em Communications on Pure and Applied Mathematics}, 13:49--54,
  1960.

\bibitem{enwiki22-NP-completeness}
{Wikipedia contributors}.
\newblock Np-completeness --- {Wikipedia}{,} the free encyclopedia.
\newblock
  \url{https://en.wikipedia.org/w/index.php?title=NP-completeness&oldid=1091328437},
  2022.
\newblock [Online; accessed 29-June-2022].

\bibitem{PTA2-xie2001regression}
Minge Xie.
\newblock Regression analysis of group testing samples.
\newblock {\em Statistics in medicine}, 20(13):1957--1969, 2001.

\bibitem{PTA3-xie2001group}
Minge Xie, Kay Tatsuoka, Jerome Sacks, and S~Stanley Young.
\newblock Group testing with blockers and synergism.
\newblock {\em Journal of the American Statistical Association},
  96(453):92--102, 2001.

\bibitem{Yao-monotonicity-1988}
F.~K. Yao, Y. C.;~Hwang.
\newblock A fundamental monotonicity in group testing.
\newblock {\em SIAM Journal on Discrete Mathematics}, 1, 05 1988.

\bibitem{yao_optimal_1990}
Y.C. Yao and F.K. Hwang.
\newblock On optimal nested group testing algorithms.
\newblock {\em Journal of Statistical Planning and Inference}, 24(2):167--175,
  February 1990.

\end{thebibliography}

\end{document}